\documentclass[12pt]{article} 

\usepackage{amsmath,amssymb}      
\numberwithin{equation}{section}
              
\newtheorem{thm}{Theorem}[section]

\newtheorem{prop}{Proposition}[section]
\newtheorem{rem}{Remark}[section]

\newcommand{\re}{\mathbb R}
\newcommand{\Nt}{\mathbb N}

\newcommand{\al}{\alpha}

\newcommand{\ep}{\varepsilon}

\def\<{\langle }

\newcommand{\supp}{\text{ supp }}

\newcommand{\qed}{\qquad\kern1pt   
   \vbox{\hrule height 0.6pt      
         \hbox{\vrule width 0.6pt 
               \vbox{\vskip 6pt}  
               \hskip 3pt
              \vrule width 1.3pt} 
         \hrule depth 1.3pt}     
   \kern1pt}

\def\v#1{\mbox{\boldmath $#1$}}

\title{The lifespan of classical solutions
of one dimensional wave equations
with semilinear terms\\
 of the spatial derivative}
\author{
Takiko Sasaki
\footnote{
Department of Mathematical Engineering, Faculty of Engineering, Musashino University,
3-3-3 Ariake, Koto-ku, Tokyo 135-8181, Japan./
Mathematical Institute, Tohoku University,
Aoba, Sendai 980-8578, Japan.
e-mail: t-sasaki@musashino-u.ac.jp.
},
Shu Takamatsu
\footnote{
Master course, Mathematical Institute,
Tohoku University,
Aoba, Sendai 980-8578, Japan.
email: shu.takamatsu.r8@dc.tohoku.ac.jp },
Hiroyuki Takamura
\footnote{Mathematical Institute,
Tohoku University,
Aoba, Sendai 980-8578, Japan.
e-mail: hiroyuki.takamura.a1@tohoku.ac.jp.}
}
\date{
\[
\begin{array}{ll}
\mbox{\footnotesize{\bf Keywords:}}
& \mbox{\footnotesize semilinear wave equation, one dimension, classical solution, lifespan}\\
\mbox{\footnotesize{\bf MSC2020:}}
& \mbox{\footnotesize primary 35L71, secondary 35B44}\\
\end{array}
\]
}

\begin{document}
\maketitle


\begin{abstract}
This paper is devoted to the lifespan estimates of small classical solutions
of the initial value problems for one dimensional wave equations
with semilinear terms of the spatial derivative of the unknown function.
It is natural that the result is same as the one for semilinear terms of the time-derivative.
But there are so many differences among their proofs.
Moreover, it is meaningful to study this problem in the sense that
it may help us to investigate its blow-up boundary in the near future.
\end{abstract}


\section{Introduction}
In this paper, we consider the initial value problems;
\begin{equation}
\label{initial_value_problem}
\left\{
\begin{array}{ll}
	 u_{tt}-u_{xx}=|u_x|^p
	&\mbox{in}\quad \re\times(0,T),\\
	u(x,0)=\ep f(x),\ u_t(x,0)=\ep g(x),
	& x\in\re,
\end{array}
\right.
\end{equation}
where $p>1$, and $T>0$.
We assume that $f$ and $g$ are given smooth functions of compact support
and a parameter $\ep>0$ is \lq\lq small enough".
We are interested in the lifespan $T(\ep)$, the maximal existence time,
of classical solutions of (\ref{initial_value_problem}).
Our result is that there exists positive constants $C_1,C_2$ independent of $\ep$ such that
$T(\ep)$ satisfies 
\begin{equation}
\label{result}
C_1\ep^{-(p-1)}\le T(\ep) \le C_2\ep^{-(p-1)}.
\end{equation}
We note that, even if $|u_x|^p$ is replaced with $|u_t|^p$, (\ref{result}) still holds.
Such a result is due to Zhou \cite{Zhou01} for the upper bound of $T(\ep)$,
and Kitamura, Morisawa and Takamura \cite{KMT23} for the lower bound of $T(\ep)$.
\par
As model equations to ensure the optimality of the general theory
for nonlinear wave equations by Li, Yu and Zhou \cite{LYZ91, LYZ92},
the nonlinear term $|u_t|^p$ is sufficient to be studied
except for the \lq\lq combined effect" case.
See Morisawa, Sasaki and Takamura \cite{MST}
and Kido, Sasaki, Takamatsu and Takamura \cite{KSTT}
for this direction with a possibility to improve the general theory.
But it is quite meaningful to deal with also $|u_x|^p$
because their proofs are technically different from each others.
Moreover, there is no result on its blow-up boundary
due to lack of the monotonicity of the solution,
while the one for $|u_t|^p$ is well-studied by Sasaki \cite{Sasaki18,Sasaki21},
and Ishiwata and Sasaki \cite{IS20a,IS20b}.
See Remark \ref{rem:non-definite} below.
It is also remarkable that it can be studied if the nonlinear term has a special form
of both $u_t$ and $u_x$. See Sasaki \cite{Sasaki22} for this direction.
Our research may help us to study the blow-up boundary for the equation
in (\ref{initial_value_problem}) near future.
\par
This paper is organized as follows. In the next section, the preliminaries are introduced.
Moreover, (\ref{result}) are divided into two theorems. 
Section 3 is devoted to the proof of the existence part, the lower bound of $T(\ep)$, of (\ref{result}). 
The main strategy is the iteration method for the system of integral equations for $(u,u_x)$,
which is essentially due to Kitamura, Morisawa and Takamura \cite{KMT23}.
They employed it for the system of integral equations for $(u,u_t)$
to construct a classical solution of the wave equation with nonlinear term $|u_t|^p$,
which is originally introduced by John \cite{John79}.
In the section 4, following Rammaha \cite{Rammaha95, Rammaha97}.
we prove the blow-up part, the upper bound of $T(\ep)$, of (\ref{result}).
We note that the method to reduced to $u$-closed integral inequality
by Zhou \cite{Zhou01} for the nonlinear term $|u_t|^p$ cannot be applicable to (\ref{initial_value_problem})
because a time delay appears in the reduced ordinary differential inequality.
Rammaha \cite{Rammaha95, Rammaha97} overcomes this difficulty
by employing weighted functionals along with the characteristic direction
in studyng two or three dimensional wave equations with nonlinear terms of spatial derivatives.


\section{Preliminaries and main results}
Throughout this paper, we assume that the initial data
$(f,g)\in C_0^2(\re)\times C^1_0(\re)$ satisfies
\begin{equation}
\label{supp_initial}
\mbox{\rm supp }f,\ \mbox{supp }g\subset\{x\in\re: |x|\le R\},\quad R\geq1. 
\end{equation}
Let $u$ be a classical solution of (\ref{initial_value_problem}) in the time interval $[0,T]$.
Then the support condition of the initial data, (\ref{supp_initial}), implies that
\begin{equation}
\label{support_sol}
\mbox{supp}\ u(x,t)\subset\{(x,t)\in\re\times[0,T]:|x|\leq t+R\}.
\end{equation}
For example, see Appendix of John \cite{John_book} for this fact.

\par
It is well-known that $u$ satisfies the integral equation;
\begin{equation}
\label{u}
u(x,t)=\ep u^0(x,t)+L(|u_x|^p)(x,t),
\end{equation}
where $u^0$ is a solution of the free wave equation with the same initial data,
\begin{equation}
\label{u^0}
u^0(x,t):=\frac{1}{2}\{f(x+t)+f(x-t)\}+\frac{1}{2}\int_{x-t}^{x+t}g(y)dy,
\end{equation}
and a linear integral operator $L$ for a function $v=v(x,t)$ in Duhamel's term is defined by
\begin{equation}
\label{nonlinear}
L(v)(x,t):=\frac{1}{2}\int_0^tds\int_{x-t+s}^{x+t-s}v(y,s)dy.
\end{equation}
Then, one can apply the time-derivative to (\ref{u}) to obtain
\begin{equation}
\label{u_t}
u_t(x,t)=\ep u_t^0(x,t)+L'(|u_x|^p)(x,t)
\end{equation}
and
\begin{equation}
\label{u^0_t}
u_t^0(x,t)=\frac{1}{2}\{f'(x+t)-f'(x-t)+g(x+t)+g(x-t)\},
\end{equation}
where $L'$ for a function $v=v(x,t)$ is defined by
\begin{equation}
\label{nonlinear_derivative}
L'(v)(x,t):=\frac{1}{2}\int_0^t\{v(x+t-s,s)+v(x-t+s,s)\}ds.
\end{equation}
Therefore, $u_t$ is expressed by $u_x$. On the other hand, applying the space-derivative to (\ref{u}),
we have
\begin{equation}
\label{u_x}
u_x(x,t)=\ep u_x^0(x,t)+\overline{L'}(|u_x|^p)(x,t)
\end{equation}
and
\begin{equation}
\label{u^0_x}
u_x^0(x,t)=\frac{1}{2}\{f'(x+t)+f'(x-t)+g(x+t)-g(x-t)\},
\end{equation}
where $ \overline{L'} $ for a function $v=v(x,t)$ is defined by
\begin{equation}
\label{nonlinear_derivative_conjugate}
\overline{L'}(v)(x,t):=
\frac{1}{2}\int_0^t\{v(x+t-s,s)-v(x-t+s,s)\}ds.
\end{equation}

\begin{rem}
\label{rem:non-definite}
In view of (\ref{u_x}), it is almost impossible to obtain a point-wise positivity of $u_x$.
This fact prevents us from studying its blow-up boundary as stated in Introduction.
\end{rem}

Moreover, applying one more time-derivative to (\ref{u_x}) yields that
\begin{equation}
\label{u_xt}
u_{xt}(x,t)=\ep u_{xt}^0(x,t)+{L'}(p|u_x|^{p-2}u_xu_{xt})(x,t)
\end{equation}
and
\begin{equation}
\label{u^0_xt}
u_{xt}^0(x,t)=\frac{1}{2}\{f''(x+t)-f''(x-t)+g'(x+t)+g'(x-t)\}.
\end{equation}
Similarly, we have that
\[
u_{tt}(x,t)=\ep u^0_{tt} + |u_x|^p(x,t) + \overline{L'}(p|u_x|^{p-2}u_xu_{xt})(x,t).
\]
Therefore, $u_{tt}$ is expressed by $u_x,u_{xt}$ and so is $u_{xx}$.

\par
First, we note the following fact.

\begin{prop}
\label{prop:system}
Assume that $(f,g)\in C_0^2(\re)\times C_0^1(\re)$.
Let $w$ be a $C^1$ solution of (\ref{u_x}) in which $u_x$ is replaced with $w$.
Then, 
\[
u(x,t):= \int^x_{-\infty}w(y,t)dy
\]
is a classical solution of (\ref{initial_value_problem}) in $\re\times[0,T]$.
\end{prop}
\par\noindent
{\bf Proof.} This is easy along with the computations above in this section. 
\hfill$\Box$

\vskip10pt
\par
Our results are divided into the following two theorems.

\begin{thm}
\label{thm:lower-bound}
Assume (\ref{supp_initial}).
Then, there exists a positive constant $\ep_1=\ep_1(f,g,p,R)>0$ such that
a classical solution $u\in C^2(\re\times[0,T])$ of (\ref{initial_value_problem}) exists
as far as $T$ satisfies
\begin{equation}
\label{lower-bound}
T \leq C_1\ep^{-(p-1)}
\end{equation}
where $0<\ep\leq\ep_1$, and $C_1$ is a positive constant independent of $\ep$.
\end{thm}

\begin{thm}
\label{thm:upper-bound}
Assume (\ref{supp_initial}) and
\begin{equation}
\label{asummption1}
f(x), g(x)\geq0,\ \mbox{and}\ f(x)\not\equiv0.
\end{equation}
Then, there exists a positive constant $\ep_2=\ep_2(f,p,R)>0$  such that
any classical solution of (\ref{initial_value_problem}) in the time interval $[0,T]$ cannot exist
as far as $T$ satisfies
\begin{equation}
\label{upperr-bound}
T>C_2\ep^{-(p-1)}
\end{equation}
where $0<\ep\leq\ep_2$, and $C_2$ is a positive constant independent of $\ep$.
\end{thm}

The proofs of above theorems are given in following sections.


\section{Proof of Theorem \ref{thm:lower-bound}}
\par
According to Proposition \ref{prop:system}, we shall construct a $C^1$ solution
of (\ref{u_x}) in which $u_x=w$ is the unknown function.
Let $\{w_j\}_{j\in \Nt}$ be a sequence of $C^1(\re\times[0,T])$ defined by
\begin{equation}
\label{w_j}
\left\{
\begin{array}{l}
w_{j+1}=\ep u_x^0+\overline{L'}(|w_j|^p), \\
w_1=\ep u_x^0.
\end{array}
\right.
\end{equation}
Then, in view of (\ref{u_xt}), $(w_j)_t$ has to satisfy
\begin{equation}
\label{w_j_t}
\left\{
\begin{array}{l}
(w_{j+1})_t=\ep u_{xt}^0+{L'}(p|w_j|^{p-2}w_j(w_j)_t), \\
(w_1)_t=\ep u_{xt}^0,
\end{array}
\right.
\end{equation}
so that the functional space in which $\{w_j\}$ converges is
\[
X:=\{w\in C^1(\re\times[0,T]): \|w\|_X<\infty, \mbox{ supp }w \subset \{(x,t)\in\re\times[0,T] : |x|\leq t+R\}\}
\]
which is equipped with a norm
\[
\|w\|_X:=\|w\|+\|w_t\|
\]
where
\[
\|w\|:= \sup_{(x,t)\in\re\times[0,T]}|w(x,t)|.
\]
We note that (\ref{u_x}) implies that
\[
\begin{array}{l}
\supp w_j\subset\{(x,t)\in\re\times[0,T] : |x|\leq t+R\}\\
\Longrightarrow\supp w_{j+1}\subset\{(x,t)\in\re\times[0,T] : |x|\leq t+R\}.
\end{array}
\]
\par
The following lemma provides us a priori estimate.
\begin{prop}
\label{prop:apriori}
Let $w\in C(\re\times[0,T])$ and supp\ $w\subset\{(x,t)\in\re\times[0,T]:|x|\leq t+R\}$. Then 
\begin{equation}
\label{apriori}
\|L'(|w|^p)\|\leq C\|w\|^p(T+R)
\end{equation}
where $C$ is a positive constant independent of $T$ and $\ep$.
\end{prop}

\par\noindent
{\bf Proof.} The proof of Proposition \ref{prop:apriori} is completely same as
the one of Proposition 3.1 in Morisawa, Sasaki and Takamura \cite{MST}.
\hfill$\Box$

\vskip5pt
\par
Let us continue to prove Theorem \ref{thm:lower-bound}.
Set
\[
M:=\sum_{\al=0}^2\|f^{(\al)}\|_{L^\infty(\re)} + \sum_{\beta=0}^1\|g^{(\beta)}\|_{L^\infty(\re)}.
\]
\vskip10pt
\par\noindent
{\bf The convergence of the sequence $\v{\{w_j\}}$.}
\par
First we note that $\|w_1\| \leq M\ep$ by (\ref{u^0_x}). (\ref{w_j}) and (\ref{apriori}) yield that
\[
\|w_{j+1} \| \leq M \ep +  C \|w_j\|^p(T+R)
\]
because it is trivial that 
\[
|\overline{L'}(v)| \leq L'(|v|).
\]
Therefore the boundedness of $\{w_j\}$, i.e. 
\begin{equation}
\label{bound_w}
\|w_j\| \leq 2M\ep \quad(j\in\Nt),
\end{equation}
follows from
\begin{equation}
\label{condi1}
C(2M\ep)^p (T+R)\leq M\ep.
\end{equation}
Assuming (\ref{condi1}), one can estimate $\|w_{j+1}-w_j\|$ as follows.
\[\begin{split}
\|w_{j+1}-w_j\| 
& =\|\overline{L'}(|w_{j}|^p-|w_{j-1}|^p ) \| \leq \|{L'}(||w_{j}|^p-|w_{j-1}|^p |) \| \\
&\leq 2^{p-1}p\|{L'}((|w_{j}|^{p-1}+|w_{j-1}|^{p-1} )|w_{j}-w_{j-1}|) \| \\
&\leq  2^{p-1}pC(\|w_{j}\|^{p-1}+\|w_{j-1}\|^{p-1})(T+R)\|w_{j}-w_{j-1}\| \\
&\leq 2^{p}pC(2M\ep)^{p-1}(T+R)\|w_{j}-w_{j-1}\|.
\end{split}\]
Therefore the convergence of $\{w_j\}$ follows from
\[
\|w_{j+1}-w_j\| \leq \frac{1}{2}\|w_j-w_{j-1}\| 
\]
provided (\ref{condi1}) and
\begin{equation}
\label{condi2}
2^{p}pC(2M\ep)^{p-1}(T+R) \leq \frac{1}{2}
\end{equation}
are fulfilled.

\vskip10pt
\par\noindent
{\bf The convergence of the sequence $\v{\{(w_j)_t\}}$.}
\par
First we note that $\|(w_1)_t\| \leq M\ep$ by (\ref{u^0_xt}).
Assume that (\ref{condi1}) and (\ref{condi2}) are fulfilled.
Since (\ref{w_j_t}) and (\ref{apriori}) yield that
\[\begin{split}
\|(w_{j+1})_t\|
&\leq M\ep+ \|{L'}(p|w_j|^{p-2}w_j(w_j)_t)\| \\
&\leq M\ep +  \|{L'}(p|w_j|^{p-1}|(w_j)_t|)\| \\
&\leq M\ep + pC\|w_j\|^{p-1}(T+R)\|(w_j)_t\| \\
&\leq M\ep + pC(2M\ep)^{p-1}(T+R)\|(w_j)_t\|, \\
\end{split}\]
the boundness of $\{(w_j)_t\}$, i.e.
\[
\|(w_j)_t\|\leq 2M\ep,
\]
follows as long as it is fulfilled that
\begin{equation}
\label{condi3}
pC(2M\ep)^{p-1} (T+R)\leq 1.
\end{equation}
Assuming (\ref{condi3}), one can estimate $\{(w_{j+1})_t-(w_j)_t\}$ as follows.
Noting that
\[
\begin{split}
&||w_j|^{p-2}w_j-|w_{j-1}|^{p-2}w_{j-1}| \\
&\leq
\left\{
\begin{array}{ll}
(p-1)2^{p-2}(|w_j|^{p-2}+|w_{j-1}|^{p-2})|w_j-w_{j-1}| & \mbox{for}\ p\ge2,\\
2|w_j-w_{j-1}|^{p-1} & \mbox{for}\ 1<p<2,
\end{array}
\right.
\end{split}
\]
we have
\[
\begin{array}{l}
\|(w_{j+1})_t-(w_j)_t\| 
= \|L'(p|w_j|^{p-2}w_j(w_j)_t -p|w_{j-1}|^{p-2}w_{j-1}(w_{j-1})_t) \| \\
\leq p \|L'(|w_j|^{p-1} | (w_j)_t -(w_{j-1})_t |) \|
+ p\|L'(|  |w_j|^{p-2}w_j -|w_{j-1}|^{p-2}w_{j-1}||(w_{j-1})_t | ) \| \\
\leq pC \|w_j\|^{p-1}(T+R)  \|(w_j)_t -(w_{j-1})_t \| \\
\quad + \left\{
\begin{array}{ll}
L'(p(p-1)2^{p-2}(|w_j|^{p-2}+|w_{j-1}|^{p-2})|w_j-w_{j-1}||(w_{j-1})_t | & \mbox{for}\ p\geq2,\\
L'(2p|w_j-w_{j-1}|^{p-1}|(w_{j-1})_t |) & \mbox{for}\ 1<p<2
\end{array}
\right. \\
\leq pC \|w_j\|^{p-1} (T+R) \|(w_j)_t -(w_{j-1})_t \| \\
\quad + \left\{
\begin{array}{ll}
p(p-1)2^{p-2}C(\|w_j\|^{p-2}+\|w_{j-1}\|^{p-2})\|w_j-w_{j-1}\|\|(w_{j-1})_t \|& \mbox{for}\ p\geq2,\\
2pC\|w_j-w_{j-1}\|^{p-1}\|(w_{j-1})_t \| & \mbox{for}\ 1<p<2
\end{array}
\right.  \\
\leq pC (2M\ep)^{p-1} (T+R) \|(w_j)_t -(w_{j-1})_t \| + O\left(\displaystyle\frac{1}{2^{j\min(p-1,1)}}\right).
\end{array}
\]
Therefore, we obtain the convergence of $\{(w_j)_t\}$ provided
\begin{equation}
\label{condi4}
pC(2M\ep)^{p-1} (T+R)\leq \frac{1}{2}.
\end{equation}

\vskip10pt
\par\noindent
{\bf Continuation of the proof.}
\par
The convergence of the sequence $\{w_j\}$ to $w$ in the closed subspace of $X$
satisfying $\|w\|, \|w_t\| \leq 2M\ep$
is established by (\ref{condi1}), (\ref{condi2}), (\ref{condi3}), and (\ref{condi4}), which follow from
\[
2^{p+1}pC(2M)^{p-1}\ep^{p-1}(T+R) \leq 1.
\]
Therefore the statement of Theorem \ref{thm:lower-bound} is established with
\[
\ep_1:=(2^{p+2}pC(2M)^{p-1}R)^{-1/(p-1)},\quad
C_1:= 2^{p+1}pC(2M)^{p-1}
\]
because $R\leq (2C_1)^{-1}\ep^{-(p-1)}$ holds for $0<\ep\leq\ep_1$.
\hfill$\Box$


\section{Proof of Theorem \ref{thm:upper-bound}}
\par
Following Rammaha \cite{Rammaha95}, set
\[
H(t):=\int^t_0(t-s)ds\int^{s+R}_{s+R_0}u(x,s)dx
\]
where $R_0$ is some fixed point with $0<R_0<R$.
We may assume that there exists a point $x_0\in(R_0,R)$ such that
$f(x_0)>0$ because of the assumption (\ref{asummption1}) and of a possible shift of $x$-variable.

\par
Then it follows that
\begin{equation}
\label{H''}
\begin{array}{ll}
H''(t)
&\displaystyle=\int^{t+R}_{t+R_0}u(x,s)dx\\
&\displaystyle=\frac{\ep}{2}\int^{t+R}_{t+R_0}
\left\{f(x+t)+f(x-t)+\frac{1}{2}\int_{x-t}^{x+t}g(y)dy\right\}dx+\frac{1}{2}F(t),
\end{array}
\end{equation}
where
\[
F(t):=\int^{t+R}_{t+R_0}dx\int_0^tds\int_{x-t+s}^{x+t-s}|u_x(y,s)|^pdy.
\]
By virtue of (\ref{asummption1}) and (\ref{H''}), we have that
\[
H''(t)\geq \frac{\ep}{2}\int^{t+R}_{t+R_0}f(x-t)dx \geq2C_f\ep
\]
where
\[
C_f:=\frac{1}{4}\int^R_{R_0}f(y)dy>0.
\]
Integrating this inequality in [$0,t$] twicely and noting that $H'(0)=H(0)=0$, we have
\begin{equation}
\label{esti_H}
H(t)\geq C_f\ep t^2.
\end{equation}

\par
On the other hand, $F(t)$ can be rewritten as
\[
F(t)=\int_0^tds\int_{t+R_0}^{t+R}dx\int_{x-t+s}^{x+t-s}|u_x(y,s)|^pdy.
\]
From now on, we assume that
\begin{equation}
\label{R_1}
t \geq R_1:=\frac{R-R_0}{2}>0.
\end{equation}
Then, inverting the order on $(y,x)$-integral,
for $0 \leq s \leq t-R_1$, we have that
\[
\begin{split}
&\int_{t+R_0}^{t+R}dx\int_{x-t+s}^{x+t-s}|u_x(y,s)|^pdy\\
&=\left(\int_{s+R_0}^{s+R}\int_{t+R_0}^{y+t-s}
+\int_{s+R}^{2t-s+R_0}\int_{t+R_0}^{t+R}
+\int_{2t-s+R_0}^{2t-s+R}\int_{y-t+s}^{t+R}\right)|u_x(y,s)|^pdxdy\\
&\geq \int_{s+R_0}^{s+R}dy\int_{t+R_0}^{y+t-s}|u_x(y,s)|^pdx.
\end{split}
\]
Similarly, for $t-R_1\leq s \leq t$, we also have that
\[
\begin{split}
&\int_{t+R_0}^{t+R}dx\int_{x-t+s}^{x+t-s}|u_x(y,s)|^pdy\\
&=\left(\int_{s+R_0}^{2t-s+R_0}\int_{t+R_0}^{y+t-s}
+\int_{2t-s+R_0}^{s+R}\int_{y-t+s}^{y+t-s}
+\int_{s+R}^{2t-s+R}\int_{y-t+s}^{t+R}\right)|u_x(y,s)|^pdxdy\\
&\geq \int_{s+R_0}^{2t-s+R_0}dy\int_{t+R_0}^{y+t-s}|u_x(y,s)|^pdx
+\int_{2t-s+R_0}^{s+R}dy\int_{y-t+s}^{y+t-s}|u_x(y,s)|^pdx.
\end{split}
\]
Hence we obtain that
\[
\begin{split}
F(t)
&\geq \int_0^{t-R_1}ds\int_{s+R_0}^{s+R}(y-s-R_0)|u_x(y,s)|^pdy\\
&\ +\int_{t-R_1}^tds\int_{s+R_0}^{2t-s+R_0}(y-s-R_0)|u_x(y,s)|^pdy\\
&\ +\int_{t-R_1}^tds\int_{2t-s+R_0}^{s+R}2(t-s)|u_x(y,s)|^pdy.
\end{split}
\]
Therefore it follows from (\ref{R_1}) and
\[
1=\frac{y-s-R_0}{y-s-R_0} \geq \frac{y-s-R_0}{R-R_0} \geq \frac{y-s-R_0}{2t}
\]
that
\[
\begin{split}
F(t)
&\geq \int_0^{t-R_1}\frac{t-s}{t}ds\int_{s+R_0}^{s+R}(y-s-R_0)|u_x(y,s)|^pdy\\
&\ +\int_{t-R_1}^t\frac{t-s}{t}ds\int_{s+R_0}^{2t-s+R_0}(y-s-R_0)|u_x(y,s)|^pdy\\
&\ +\int_{t-R_1}^t2(t-s)ds\int_{2t-s+R_0}^{s+R}\frac{y-s-R_0}{2t}|u_x(y,s)|^pdy\\
&= \frac{1}{t}\int_0^t(t-s)ds\int_{s+R_0}^{s+R}(y-s-R_0)|u_x(y,s)|^pdy.\\
\end{split}
\]
In this way, (\ref{asummption1}), (\ref{H''}) and the estimate of $F(t)$ above yield that
\[
H''(t)\geq \frac{1}{2}F(t)\geq \frac{1}{2t}\int_0^t(t-s)ds\int_{s+R_0}^{s+R}(y-s-R_0)|u_x(y,s)|^pdy
\quad\mbox{for}\ t\geq R_1.
\]
Moreover, it follows from (\ref{support_sol}), integration by parts and H\"older's inequality that
\[
\begin{split}
|H(t)|
&=\left|\int_0^t(t-s)ds\int_{s+R_0}^{s+R}\partial_y(y-s-R_0)u(y,s)dy\right|\\
&=\left|\int_0^t(t-s)ds\int_{s+R_0}^{s+R}(y-s-R_0)u_x(y,s)dy\right|\\
&\leq \int_0^t(t-s)ds\int_{s+R_0}^{s+R}(y-s-R_0)|u_x(y,s)|dy\\
&\leq \left(\int_0^t(t-s)ds\int_{s+R_0}^{s+R}(y-s-R_0)|u_x(y,s)|^p dy\right)^{1/p}I(t)^{1-1/p},\\
\end{split}
\]
where
\[
I(t):=\int_0^t(t-s)ds\int_{s+R_0}^{s+R}(y-s-R_0)dy=\frac{1}{4}t^2(R-R_0)^2=t^2R_1^2.
\]

\par
Hence we obtain that
\begin{equation}
\label{esti_H''}
H''(t)\geq \frac{1}{2}R_1^{-2(p-1)}t^{1-2p}|H(t)|^p\quad\mbox{for}\ t\geq R_1.
\end{equation}
Therefore, the argument in  Rammaha \cite{Rammaha95} can be applied to
(\ref{esti_H}) and (\ref{esti_H''}) to ensure that
 there exist positive constants $\ep_2=\ep_2(f,p,R)$ and $C_2$
 independent of $\ep$ such that a contradiction appears provided
\[
T>C_2\ep^{-(p-1)}
\]
holds for $0<\ep \leq \ep_2$.
The proof is now completed.
\hfill$\Box$

\section*{Acknowledgement}
\par
The first author is partially supported
by the Grant-in-Aid for Young Scientists (No. 18K13447), 
Japan Society for the Promotion of Science.
The third author is partially supported
by the Grant-in-Aid for Scientific Research (A) (No. 22H00097), 
Japan Society for the Promotion of Science.


\bibliographystyle{plain}

\end{document}